\documentclass[letterpaper, 10 pt, conference]{ieeeconf}  

\IEEEoverridecommandlockouts                              

\usepackage{graphicx} 
\usepackage{amsmath} 
\usepackage{amssymb}  

\usepackage{amsthm}

\usepackage{soul}
\usepackage[dvipsnames]{xcolor}

\theoremstyle{definition}
\newtheorem{definition}{Definition}

\newtheorem{lemma}{Lemma}

\newtheorem{problem}{Problem}
\newtheorem{remark}{Remark}

\newtheorem{assumption}{Assumption}

\newtheorem{example}{Example}

\usepackage{subcaption}

\newtheorem{objective}{Objective}

\usepackage{caption}
\DeclareCaptionType{equ}[][]

\graphicspath{{figures/}} 

\title{\LARGE \bf
Recursive Feasibility Guided Optimal Parameter Adaptation of Differential Convex Optimization Policies for Safety-Critical Systems
}

\author{Hardik Parwana$^{1}$ and Dimitra Panagou$^{2}$ \\ 
\textit{NOTE: This paper has been accepted at ICRA 2022}
\thanks{$^{1}$Hardik Parwana is with Robotics Institute, University of Michigan, Ann Arbor, MI 48109, USA
        {\tt\small hardiksp@umich.edu}}%
\thanks{$^{2}$Dimitra Panagou with the Department of Aerospace Engineering and Robotics Institute, University of Michigan, Ann Arbor, MI 48109, USA
        {\tt\small dpanagou@umich.edu}}%
}

\newcommand{\reals}{\mathbb{R}}
\newcommand{\s}{\mathcal{S}}
\newcommand{\X}{\mathcal{X}}

\newcommand{\K}{\mathcal{K}}
\newcommand{\classK}{class-$\K$ }

\newcommand{\classKinf}{\mbox{class-$\K_\infty$}}
\newcommand{\Int}{\text{Int} }

\newcommand{\eqn}[1]{\begin{align}#1\end{align}}

\begin{document}

\maketitle
\thispagestyle{empty}
\pagestyle{empty}

\begin{abstract}
Quadratic Program(QP) based state-feedback controllers, whose inequality constraints bound the rate of change of control barrier(CBFs) and lyapunov function with a class-$\mathcal{K}$ function of their values, are sensitive to the parameters of these class-$\mathcal{K}$ functions. The construction of valid CBFs, however, is not straightforward, and for arbitrarily chosen parameters of the QP, the system trajectories may enter states at which the QP either eventually becomes infeasible, or may not achieve desired performance. In this work, we pose the control synthesis problem as a differential policy whose parameters are optimized for performance over a time horizon at high level, thus resulting in a bi-level optimization routine. In the absence of knowledge of the set of feasible parameters, we develop a Recursive Feasibility Guided Gradient Descent approach for updating the parameters of QP so that the new solution performs at least as well as previous solution. By considering the dynamical system as a directed graph over time, this work presents a novel way of optimizing performance of a QP controller over a time horizon for multiple CBFs by (1) using the gradient of its solution with respect to its parameters by employing sensitivity analysis, and (2) backpropagating these as well as system dynamics gradients to update parameters while maintaining feasibility of QPs.

\end{abstract}

\section{INTRODUCTION}

Autonomous systems are expected to perform in a constrained environment and complete their tasks. The decision making is usually broken down into stability and safety objectives which can be accounted for in an optimal control formulation. Techniques exist to solve these problems, but the curse of dimensionality for Hamilton-Jacobi based methods\cite{bansal2017hamilton} or computational complexity of nonlinear optimization for Model Predictive Control type of approaches have been the bottleneck for their use in real-time implementations. As such, state feedback controllers that depend only on the current state have prevailed. 

Control Lyapunov (CLFs) functions and control barrier functions (CBFs)\cite{ames2016control}\cite{ames2014control} are two popular methods to encode stability and safety. Utilizing the dynamics of system, they restrict the rate of change of lyapunov and barrier functions by a parametric class-$\mathcal{K}$ function of their values. Quadratic Programs(QPs) can then be formulated to obtain minimum norm solution satisfying CBF condition strictly and CLF condition as best as possible\cite{ames2016control}.
However, most of the work focuses on a single barrier constraint and it is not straightforward to extend these results to multiple barrier functions. This is because it is not guaranteed that the system trajectory will reach states that will maintain feasibility of QP. This is a direct consequence of poorly chosen parameters of class-$\mathcal{K}$ function, as well as the \textit{myopic} nature of state feedback controllers. Moreover, as noted in \cite{cohen2020approximate}, these myopic controllers may also give non-optimal solutions in long term, thus achieving sub-par performance.

Note that there exist methods to combine all barriers into a single barrier 
function\cite{glotfelter2017nonsmooth, stipanovic2012monotone, panagou2013multi}. However, we argue that combination of barriers is not so intuitive to tune for, especially when they represent different physical quantities such as position and angles, whose domain and rates of change have different scales with respect to each other. 

This work addresses two issues. First, we start with parameterized QPs as policies that need to optimized. The sensitivity analysis of QPs is exploited to extract gradients of the control input with respect to its parameters. Then we consider performance over a time horizon as the criterion for adapting the parameters in presence of multiple conflicting constraints. The performance here is defined in terms of an objective function for the horizon that needs to be maximized. Therefore, this adaptation explicitly creates a link between static state-feedback controllers and their long term behavior. The resulting process is realized as a directed computational graph whose nodes represent states or functions, and edges represent the arguments to the node's functions. Parameter update is achieved by backpropagating gradients through the graph and applying constrained Gradient Descent(GD) to improve performance.
This is the first time to best of authors' knowledge that sensitivity analysis is applied to treat QP based safety-critical control design as differentiable policy that is improved over a time horizon by rolling out trajectory to future.

Secondly, the potential infeasibility of QPs at some point in future, which we call lack of recursive feasibility motivated by MPC literature, is addressed by changing parameters in direction that reduces the infeasibility margin of the constraints in QP that fail to be satisfied. Note that feasibility may be lost even in absence of control input bounds as there might not exist a control input that can satisfy all constraints simultaneously at that state.

The above methodology is made possible by advances made in several works\cite{agrawal2019differentiating}\cite{barratt2018differentiability} that study differentiating through conic optimization problems and provide conditions when this is possible. These are based on taking matrix differentials of Karush-Kuhn-Tucker(KKT)\cite{boyd2004convex} conditions of the optimization problem at its solutions. The idea to use optimization problems as differentiable layers has been realized in some recent works\cite{amos2017optnet} where instead of using standard units like ReLU as activation functions of neurons in neural networks(NNs), QPs are used to enable richer interactions between NN layers and encode any additional constraints. However, using them as a differentiable unit in dynamical system design has not been explored yet.

Compared to the above literature, this paper considers an additional issue that the optimization problem may be infeasible for some states and parameter values. Therefore unlike the non-optimization based models, an output may not even exist. However, the advantage is that just having a valid output ensures safety of the system which is unlike other policies, such as standard neural networks, where safety constraint is on output of policy. This fundamental distinction between them leads us to develop a gradient descent variation suited to QPs based on Sequential Quadratic Programming(SQP) and Feasible SQP\cite{tits2009feasible}(FSQP) literature. Earlier relevant work has appeared in \cite{xiao2020feasibility} where the proposed approach does update parameters of QP but performs an offline analysis to construct a set approximation of all feasible parameters. In our approach, such a set is not constructed as it might be expensive to compute and dependent on reference trajectories. \cite{xiao2020feasibility} also does not use policy gradients similar to our approach. We define a performance objective that represents the desired behaviour over a time horizon. We then present Recursive Feasibility Guided Gradient Descent for updating parameters of optimization problem to (1) ensure that the new parameters achieve better high-level performance than previous parameters if both are feasible over the same time horizon, or (2) if QP becomes infeasible at some point in time, then the updated parameters should keep the QP feasible for longer period of time.






\section{Preliminaries}

\subsection{Notations}
We denote the set of real numbers as $\reals$ and the non-negative real numbers as $\reals^+$. Given $x\in \reals$ and $y\in \reals^n$, $|x|$ denotes the absolute value of $x$ and $||y||$ denotes the $L_2$ norm of $y$. $\langle a,b \rangle = a^Tb$ represents the inner product between $a,b \in \reals^n$. The interior and boundary of a set $\mathcal{C}$ is denoted by $\Int(C)$ and $\partial C$. A continuous function $\alpha: [0,a)\rightarrow [0,\infty)$ is a \classK function if it is strictly increasing and $\alpha(0)=0$. Furthermore, if $a=\infty$ and $\lim_{r \rightarrow \infty} \alpha(r) = \infty$, then it is called \classKinf. Both $\frac{\partial }{\partial x}$ and $\nabla_x$ denote gradient and will be used interchangeably depending on the complexity of expressions for easy understanding.

\subsection{System Description}

Consider a control affine nonlinear dynamical system with state $x \in \mathcal{X}\subset \reals^n $, control input $u\in \mathcal{U} \subset \reals^m$ given in discrete time as
\eqn{
    x_{t+1} = f(x_t) + g(x_t)u_t \label{eq::dynamics}
}
where $f(x):\mathcal{X}\rightarrow \reals^n,g(x):\mathcal{U}\rightarrow \reals^{n\times m}$ are Lipschitz continuous and $t$ is the time index. The safety of the system is specified in terms of, possibly time varying but intersecting, safety sets $\s_i(t), i\in {1,2,..,N}$ that encode the allowable states of the system and are defined as 0-superlevel set of a smooth function $h:\mathcal{X} \times \reals^+ \rightarrow \reals$ as follows
\eqn{
        \s_i(t) & \triangleq \{ x \in \X : h_i(t,x) \geq 0 \} \\
        \partial \s_i(t) & \triangleq \{ x\in \X: h_i(t,x)=0 \} \\
        \Int (\s_i)(t) & \triangleq \{ x \in \X: h_i(t,x)>0  \}
}

The set $\mathcal{S}_i(t)$ is forward invariant if and only if $\dot{h}_i(t,x)\geq 0 \; \forall x \in \partial \mathcal{S}_i(t)$\cite{blanchini1999set} in continuous-time systems. With slight abuse of notation, we will use $\s=\cap S_i$. If $h_i$ is a zeroing-barrier function for a $\classKinf$ function $\nu_i$, then the following condition in discrete-time on $h_i$, called CBF condition, is sufficient for the  invariance of the set $\mathcal S_i$:
\eqn{
    h_i(t+1,x_{t+1}) - h_i(t,x_t) \geq \nu_i( h_i(t,x_t) ) \label{eq:CBF_cond}
}
For simplicity, we will be using $\nu_i(x)=\alpha_i x, \alpha_i\in \reals^+$ in remainder of paper but any other parametric function can also be used. 
\begin{problem}
    \label{problem::1}
Given the dynamical system (\ref{eq::dynamics})  and initial state $x_0 \in \mathcal{S}(0)$, design a controller $u(x,t)$ such that $\forall t>0$, the closed loop trajectory satisfies $x(t) \in \mathcal{S}(t)$ and $|| x_t - x^d_t||\rightarrow 0$ as $t \rightarrow \infty$.
\end{problem}
In order to achieve this objective, the following CBF-CLF-QP based approach, which we shall refer to as policy $\pi_{QP}(x_t;{\alpha_0,\alpha_1,...\alpha_N})$, is commonly used:
\begin{subequations}
        \begin{align}
            \min_{u_t,\delta} \quad & J(u) = (u_t-u_{d_t})^TP(u_t-u_{d_t}) + Q\delta^2\\
            \textrm{s.t.} \quad &  V(t+1,x_{t+1}) \leq (1-\alpha_0)V(t,x) + \delta \\
                    &       h_i(t+1,x_{t+1}) \geq (1-\alpha_i)h_i(t,x_t) ~ \forall i \in \{1,2,..,N\}
        \end{align}
        \label{eq::CBF_MDP}
\end{subequations}
Here $\alpha_i, i\in {0,1,...,N}$ are parameters of the policy, $u_{d_t}$ is the desired or nominal control input, $P\in \reals^{m\times m}, Q\in \reals^+$ are the weight factors, $V(x)$ is a given Lyapunov function, and $\delta\in \reals$ is slack variable used to relax CLF condition. In the subsequent, we will refer to all parameters together in a vector $\theta = [\alpha_0, ...\alpha_N]^T$. Note that $u_{d_t}, P$, and $Q$ may also be treated as parameters if desired so.

\begin{remark}
Refer to \cite{agrawal2017discrete},\cite{breeden2021control} for a thorough discussion of discrete time CLF and CBF conditions and their variants. We use a simple extension of their work that ensures constraint satisfaction only at sampling times. Any other variation with stronger inter-sample guarantees may be used too. 
\end{remark}
\begin{example}
\label{example::unicycle}
A follower modeled as a unicycle is supposed to follow a leader and maintain it within its Field-of-View(FoV). The Lyapunov function is defined on desired radial distance only, while three barrier functions $h_1,h_2,h_3$ encode minimum distance, maximum distance, and maximum angle w.r.t camera axis constraints. A control policy can be synthesized as
\eqn{
    \begin{aligned}
    \label{eq::unicycle_QP}
            \min \quad & u^T u + Q \delta^2 \\
            \textrm{s.t.} \quad &  V(t+1,x_{t+1}) \leq -(1-k)V(t,x_t) + \delta  \\
                          &   h_1(t+1,x_{t+1}) \geq (1-\alpha_1) h_1(t,x_t)   \\
                          &   h_2(t+1,x_{t+1}) \geq (1-\alpha_2) h_2(t,x_t)   \\
                          &   h_3(t+1,x_{t+1}) \geq (1-\alpha_3) h_3(t,x_t)   \\
    \end{aligned}
}
Here, $k,\alpha_1,\alpha_2,\alpha_3$ are parameters of the problem.
 \end{example}

Most implementations choose the parameters $\theta$ based on manual tuning and assume that $h_i$ are valid CBF that can be satisfied simultaneously at all states. However, this assumption is not true in general, and system trajectories may violate the above constraints even in presence of unbounded control inputs. Even when constraints do remain feasible, they may not allow desirable performance and there exists no automatic tuning procedure for $\theta$. The following example shows how even in simple scenarios, having multiple constraints can lead to infeasibility in future.

\begin{figure}[th!]
    
        \centering
        \includegraphics[width=0.95\linewidth]{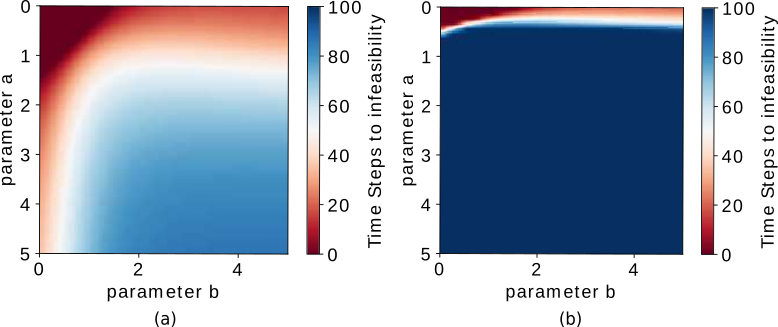}
        \caption{Time to Infeasibility: X and Y axis represent parameter values(normalized by sampling time $\Delta t$) for which trajectory is simulated. The color gradation represents the time after which problem becomes infeasible. A value of 100 implies it remains feasible for maximum possible time and value of 0 implies it is infeasible from the start. (a)c=0.3, (b)c=0.7 in Example \ref{example:car}.}
    \label{fig:1d_example}
\end{figure}

\begin{example}
\label{example:car}
Consider the motion of an autonomous car modeled as one-dimensional integrator system $\dot{x}=u$ with state $x$ and velocity input $u$. Now consider a common road scenario in which the car is caught between two vehicles. The vehicle in front is moving slower than the vehicle at back and therefore the safe set, described as no collision zone between the vehicles, is shrinking with time and vanishes in future. The objective is to stay safe for maximum possible time. Let the shrinking set be given by the constraint $x\geq t$, and $x\leq 1+c t$ where $t$ is time and $c<1$ is a constant. We can encode these two constraints with two barrier functions $h_1 = x - t$, $h_2 = 1 + c t-x$ and solve the following optimization problem
\eqn{
    \begin{aligned}
            \max \quad & u  \\
            \textrm{s.t.} \quad & h_1(t+1,x_{t+1}) \geq (1-a) h_1(t,x_t)   \\
                          &  h_2(t+1,x_{t+1}) \geq (1-b) h_2(t,x_t) 
            \label{eq::car_LP}
    \end{aligned}
}
Here, $a,b\in \reals $ are parameters. Linear Programs(LP) satisfy the same assumptions introduced for QPs in Section \ref{section::cvx_der} and therefore are still amenable to our approach. Note that eventually the safe set completely vanishes here but the response of the system before that critical time is still relevant to real scenarios. The closed loop system for different values of parameter is simulated. Fig.\ref{fig:example_path} shows the resulting trajectories and Fig.\ref{fig:1d_example} shows how time to infeasibility changes with parameter choice for two different set shrinkage rates $c$. Both figures are asymmetric w.r.t $a,b$ which is a consequence of difference in rate of change of both constraints. Figs.\ref{fig:example_path},\ref{fig:1d_example} also show why it is important to choose parameter proactively as a wrong parameter choice may lead to infeasibility and possible crash of the system.

\begin{figure}[th!]
        \centering
        \includegraphics[width=0.7\linewidth]{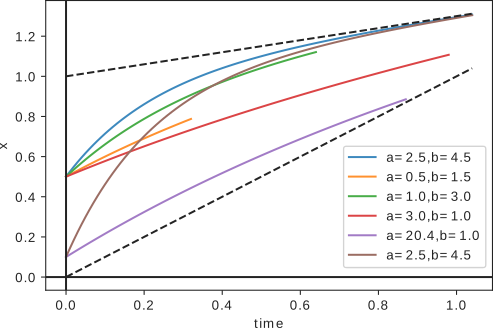}
        \caption{Plot of trajectories for different values of parameter. Trajectories end at point of infeasibility. The values of $a$ and $b$ shown are normalized by sampling time $\Delta t$. 
        }
    \label{fig:example_path}
\end{figure}
\end{example}

\subsection{Differentiable Convex Programs}
\label{section::cvx_der}

Finally for sake of completion, we mention this important result that allows us to differentiate solution of QPs with respect to its parameters. This validates our use of gradients in subsequent sections. We first make following assumptions concerning KKT conditions.
\begin{assumption}
\label{assump::slater}
(Strong Duality) Slater's condition holds for (\ref{eq::CBF_MDP})
\end{assumption}

\begin{assumption}
\label{assump::differentiable}
Let $e(t,u_{t},x_{t},\theta)\in \reals^{N+1}\leq 0$ represent the vector of all inequality constraints in (\ref{eq::CBF_MDP}). The functions $e_i$ and the objective $J(\pi)$ are twice differentiable in $u_{t}$, and $e(u_t,x_t,\theta), \nabla_u e(u_t,x_t,\theta)$ are continuously differentiable in $(x_t,\theta)$. Note that for QP, $(x_t,\theta)$ are parameters and $u_t$ is a variable. This assumption would usually be satisfied for most systems of interest that have differential dynamics and a smooth objective function.
\end{assumption}
\begin{assumption}
\label{assump::complementarity}
Strict complementarity holds for (\ref{eq::CBF_MDP}).
\end{assumption}

\begin{remark}
Note that in many control problems, the solution of QP may lie exactly on the boundary of constraint inequalities and Assumptions \ref{assump::slater},\ref{assump::complementarity} may not hold and gradient may be undefined. However even in these cases, relaxed KKT conditions\cite{biegler2010nonlinear} can be used to get gradient in the limiting case\cite{gros2019towards}. Alternatively, a more thorough analysis based on directional differential along a parameter trajectory may be possible but is out of scope of this paper. A simple way to ensure differentiability of QP w.r.t its parameters is to ensure that $u_t = u_{d_t}$ does not make the constraints active.
\end{remark}

\begin{lemma}
\textbf{(\cite{agrawal2019differentiating},\cite{barratt2018differentiability}Theorem3.1)}If Assumptions \ref{assump::slater},\ref{assump::differentiable}, and \ref{assump::complementarity} hold, then the QP solution is a locally single valued function around its solution $u_t$ and thus is continuously differentiable in the neighborhood of ($x_t,\theta$).(Refer to \cite{barratt2018differentiability} for a more insightful condition).
\end{lemma}

\section{Problem Statement}
A typical objective to be achieved over a time horizon is to maximize the summation of stage-wise cost, which we call the \textit{Reward} $r$, subject to constraints (\ref{eq:CBF_cond}b,\ref{eq:CBF_cond}c).This second optimization routine at the high-level is given as follows:
\begin{equation}
        \begin{aligned}
            \max_{\theta} \quad & R(x_1,\theta) = \sum_t^T  r(x_t,u_t) \\
            \textrm{s.t.} \quad & u_t = \pi_{QP}(x_t;\theta) \\
             & x_{t+1} = f(x_t) + g(x_t)u_t \\
             & \forall t \in \{1,2,..,T\}
        \end{aligned}
        \label{eq::high_level}
\end{equation}
The QP controller is said to be implemented at \textit{low-level}, and the high-level optimizes parameters that are used by low-level.  $r$ can encode desired stability and/or safety objectives.
The objective is to design an update rule for $\theta$ without losing feasibility of QP over period of interest, i.e., ensuring that a trajectory exists over the same horizon for updated $\theta$. 

\begin{definition}
A parameter $\theta$ is called feasible if starting at state $x$, a solution to QP exists $\forall t$. 
\end{definition}
\begin{definition}
A parameter $\theta$ is called T-feasible at state $x$ if a solution to QP exists $\forall t \leq T$ and no solution exists for $t=T+1$.
\end{definition}
In practice, it is sufficient to call a parameter feasible if the system remains feasible for the period of interest. Depending on whether we have feasibility or T-feasibility, we now specify two objectives to be solved in this paper:
\begin{objective}
  Given a feasible $\theta$ design an update rule $\theta^+=F(\theta)$ such that $\theta^+$ is feasible and $R(x_t,\theta^+)\geq R(x_t,\theta)$.
\end{objective}

\begin{objective}
    Given a T-feasible $\theta$ design an update rule $\theta^+=F(\theta)$ such that $\theta^+$ is (T+1)-feasible.
\end{objective}



\section{Methodology}
Firstly, we design an update rule for improving performance of a feasible parameter. This is done by computing gradient of the objective $\nabla_\theta J(x,\theta)$, projecting it to the set of feasible directions of $\theta$, and then updating in the resulting direction with a learning rate $\beta$. Here, the set of feasible directions is computed such that a change in this direction will maintain feasibility of the QPs and guarantee existence of trajectory at all times (or upto a large horizon in practice). This is done in manner similar to FSQP where parameters are updated so that each iterate of gradient descent, not just the last one, is a feasible solution.

Secondly, we design a rule to update a T-feasible parameter to (T+1)-feasible parameter. This is done by first identifying the constraint at T+1, 
that needs to be relaxed the least of all other constraints so that the QP at T+1 becomes feasible. A direction that reduces the infeasibility margin of this constraint is then computed 
with respect to $\theta$. This direction is then projected to set of feasible directions of QPs from t=1 to t=T and the resulting direction is used for update with learning rate $\beta$. 
The next two subsections discuss how we can compute gradients by backpropagating and then how we can use them in FSQP framework to update parameters.

\subsection{Gradient through Backpropagation}
\label{section::backprop}
\begin{figure}[htp]
    \centering
    \includegraphics[scale=0.3]{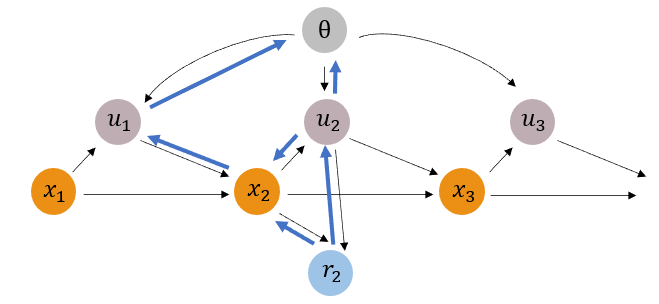}
    \caption{Nodes represent functions, specified here by dynamics or QP controller, and edges represent arguments to the function. Thick blue arrows show the flow of information when computing gradient of second reward $r_2=R(x_2,u_2)$ with respect to $\theta$ through backpropagation.}
    \label{fig:graph}
\end{figure}

The system can be visualized as a computational graph(see Fig.\ref{fig:graph}) with nodes representing independent variables and functions, and incoming edges representing arguments to the function and the flow of information. The gradient of each edge can be obtained either from the system dynamics or through the sensitivity analysis of QP. Since the dynamics is deterministic and known, we can roll-out trajectory to future and compute the gradient of objective function with respect to $\theta$ and initial state $x_1$. As an example, for $T=2$,
\eqn{
   \frac{\partial R(x_1,\theta)}{\partial \theta} = \frac{\partial r(x_1,u_1)}{\partial \theta} + \frac{\partial r(x_2,u_2)}{\partial \theta}
}
Here, we have
\eqn{
   \frac{\partial r(x_2,u_2)}{\partial \theta} = \underbrace{\frac{\partial r_2 }{\partial x_2}}_{known}\frac{\partial x_2}{\partial \theta} + \underbrace{\frac{\partial r_2}{\partial u_2}}_{known}\frac{\partial u_2}{\partial \theta}
}
Since $x_2 = f(x_1) + g(x_1)u_1$ is a function of $x_1$ and $u_1$
\eqn{
  \frac{\partial x_2}{\partial \theta} = \underbrace{\frac{\partial x_2}{\partial x_1 }}_{\textrm{from dynamics}}\underbrace{\frac{\partial x_1 }{\partial \theta}}_{\textrm{$x_1$ constant}} + \underbrace{\frac{\partial x_2}{\partial u_1}}_{\textrm{from dynamics}}\frac{\partial u_1}{\partial \theta}
}
$\frac{\partial u_1}{\partial \theta}, \frac{\partial u_2}{\partial \theta}$ are gradients of QP's solution w.r.t the parameters which can be obtained using available libraries\cite{diamond2016cvxpy}. 
This recursive analysis is amenable to backpropagation type of derivative computation and can easily be implemented in libraries., like pytorch and tensorflow, capable of constructing computational graphs. The gradient of other quantities of interest, such as barrier function, can be obtained similarly.


\vspace{-5mm}
\subsection{Recursive Feasibility Guided Gradient Descent(RFGGD)}

We provide two different algorithms for the two objectives.

\subsubsection{Case 1: $\theta_t$ is feasible}
\label{algo::case1}
In this case, we are only concerned with optimizing the high level objective $R(x_1,\theta)$ 
An unconstrained update with $\theta_{t+1} = \theta_t + \beta \nabla_\theta R$ would not guarantee that the QPs will be feasible for $\theta_{t+1}$. Also, when we update $\theta$, the resulting trajectory will be different from previous $\theta$ over the time horizon and therefore we must design an update rule that guarantees a feasible trajectory. This requires us to consider variations of
all constraints across the horizon together. 
Below, similar to SQP, we show how to get feasible directions upto first-order approximations of current trajectory, but any FSQP result\cite{tits2009feasible} can be used to augment this step to give feasibility guarantees in general.

The high level optimization (\ref{eq::high_level}) is a nonlinear optimization problem as the barrier and lyapunov constraints are not linear in state of the system along a horizon. Let $e(t,u_t,x_t,\theta) \in \reals^{N+1} \geq 0$ represent all inequality constraints at time $t$. We then use first order approximation of the effect of changing parameter on this constraint by linearizing it about current solution $x_t,u_t$. The following is then enforced
\eqn{
    e(t,u_t,x_t;\theta) +  \nabla_\theta e(t,u_t,x_t;\theta) d_\theta \geq 0
}
where $d_\theta$ is a possible perturbation in $\theta$ such that above equation is satisfied. Here, the gradient $\nabla_\theta e$ is computed by backpropagation as $x_t,u_t$ are also variables that recursively depends only on initial state $x_1$ and $\theta$. 
The set of feasible directions, upto first order approximation, is thus given by
\begin{align*}
    \mathcal{F}(x_1,\theta) &= \{d_{\theta}~ | ~ e(t,u_t,x_t;\theta) + \nabla e(t,u_t,x_t;\theta) d_{\theta} \geq 0, \\
    & \quad ~~~~~~~~~~~~~~~~~~~~~~~~~~~~~  \forall t\in \{1,2,..,T \} ~ \}
\end{align*}

Finally, we can find an update direction for Objective 1 by solving the following LP
\begin{equation}
    \begin{aligned}
        \max_{d_\theta} \quad & \langle \nabla_\theta R(x_1,\theta), d_\theta \rangle \\
        &d_\theta \in \mathcal{F}(x_1,\theta)
    \end{aligned}
    \label{eq::feasible_update}
\end{equation}

where the objective corresponds to maximizing the projection of ascent direction of the high level objective $\nabla_\theta R(x_1,\theta)$ along the feasible direction $d_\theta$ of parameter update. A solution to this problem always exists as $d_\theta=0$ satisfies the inequality. Once a suitable direction has been found, parameters can be updated as follows
\eqn{
\theta_{t+1} = \theta_t + d_\theta
\label{eq::projected_GD}
}
Note again that this result guarantees feasibility maintenance only up to first order approximation. More precise methods from FSQP may be used here or a learning rate $\beta$ might be used in Eq.(\ref{eq::projected_GD}) to promote small changes only.

\subsubsection{Case 2: $\theta_t$ is T-feasible}
\label{algo::case2}
The objective here is to make the QP at $T+1$ feasible. We first find the constraint at $T+1$ that needs to be relaxed by the least amount to make the QP feasible. There are many ways to do so and a simple one would be to add slack variables to all constraints and find a feasible solution minimizing the sum of squares of slacks. The non-zero slack values would then correspond to the limiting constraints, with index belonging to set $L\subset \{1,2,..,N+1\}$, denoted as $e_j(t=T,x_T,u_T,\theta)\geq 0, j\in L$, which cannot be satisfied with the current parameter $\theta$. Then we would like to update $\theta$ so that $e_j$ can be increased. Suppose $e_j = C_1(t=T,x_T,\theta) + C_2^T(t=T,x_T,\theta)u_T$, for some suitable functions $C_1,C_2$, which is always the case for CBFs with control affine systems. Then we find the ascent direction $d_j$ of $C_1$ so that $e_j$ can potentially be raised to $>0$. This is again done with backpropagation as $x_T$ is recursively related to $x_1$ and $\theta$. The update direction for Objective 2 is obtained by projecting $d_j$ to the set of feasible direction of $\theta$
which is computed in same manner as before. 
\begin{equation}
    \begin{aligned}
        \max_{d_\theta} \quad & \langle d_j, d_\theta \rangle \\
        &d_\theta \in \mathcal{F}(x_1,\theta)
    \end{aligned}
    \label{eq::infeasible_update}
\end{equation}
\begin{remark}
In practice, for Case 1, a single update of GD can be done at each time step. For Case 2, we can perform multiple GD updates until $T+1$ becomes feasible. Since each update is based on single trajectory roll-out and backpropagation, they can be performed very fast.
\end{remark}
\begin{remark}
The dynamics equation is not explicitly considered in computation of feasibility directions as (1)CBFs have conveniently combined the dynamics with state constraint in a single inequality and (2) the backpropagation does take into account dynamics while traversing the computational graph.
\end{remark}

\section{Simulation Results}
We present results for two case studies. The first one corresponds to the autonomous car model of Example \ref{example:car}, and applies our algorithm to increase the horizon over which the system is feasible. The second one corresponds to the leader-follower problem of Example \ref{example::unicycle}, and shows how parameter adaptation can improve the high-level objective.
The optimization problems are solved in python with cvxpy\cite{diamond2016cvxpy} interface and SCS solver. cvxpy also returns the gradient of QP solution with respect to its parameters. The dynamics of robotic agent and cvxpy optimization are made part of PyTorch layers and once the computation graph is made, backpropagation is done with PyTorch's Autograd feature\footnote{The videos and the code implementing the algorithm can be found at https://github.com/hardikparwana/Safe-Learning-DASC}. 

\subsection{Autonomous car with shrinking safe set}
\begin{figure}[th!]
        \centering
        \includegraphics[width=1.0\linewidth]{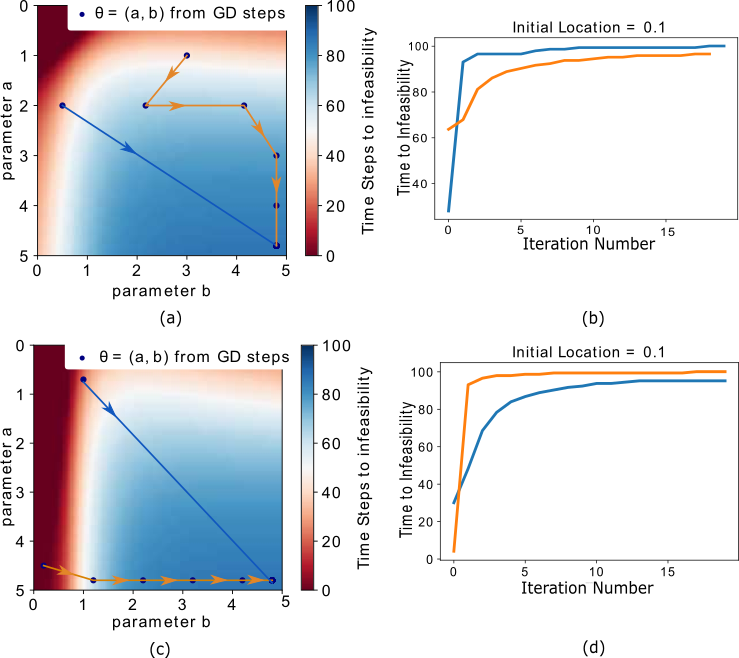}
        \caption{GD iterates: All figures are for c=0.3 in Example \ref{example:car}. Initial state is 0.5 for (a),(b), and 0.1 for (c),(d). Each shows results for two different initialization of $\theta$. (a), (c) show the change in $\theta$. Points on same GD process are connected by line and direction of update marked with arrows. (b), (d) show plots of time to infeasibility with GD iteration number. Parameter values have been clipped to range (0,5) to better visualize infeasible parameters near top and left edges of the box. The values of $a$ and $b$ shown are normalized by $\Delta t$.}
    \label{fig:car_plots}
\end{figure}

Fig.\ref{fig:car_plots} shows the results from applying algorithm from Section \ref{algo::case2} to autonomous car problem. The constrained GD is able to update parameter while ensuring that the system remains feasible for the same or longer time horizon.

\subsection{Leader-Follower Problem}
The position and velocity of the leader is assumed to known to the follower. The FoV constraint specifies a time varying safe set obtained as intersection of following three barrier conditions $h_1 = s^2 - s_{min}^2\geq 0, h_2 = s_{max}^2 - s^2 \geq 0$, and $h_3 = b - \cos(\gamma) \geq 0$ where $s$ is the distance between the follower and the leader, $b$ is the bearing vector from the follower to the leader, and $\gamma$ is the FoV angle of the follower's camera. The parameters to adapt are as mentioned in Example \ref{example::unicycle}. Additionally, a CLF condition with lyapunov function $V=(s-s_d)^2$, $s_d$ being desired distance between follower and leader, is also added to promote convergence of leader to the center of FoV. The QP controller to be solved is given by Eq.(\ref{eq::unicycle_QP}). The control inputs $u$ of follower are its linear and angular velocity. The reward function $r$ is chosen to be the smooth minimum of the three barrier functions and is higher if the leader is maintained at center of the FoV. 

\begin{figure}[th!]
        \centering
        \includegraphics[width=1.0\linewidth]{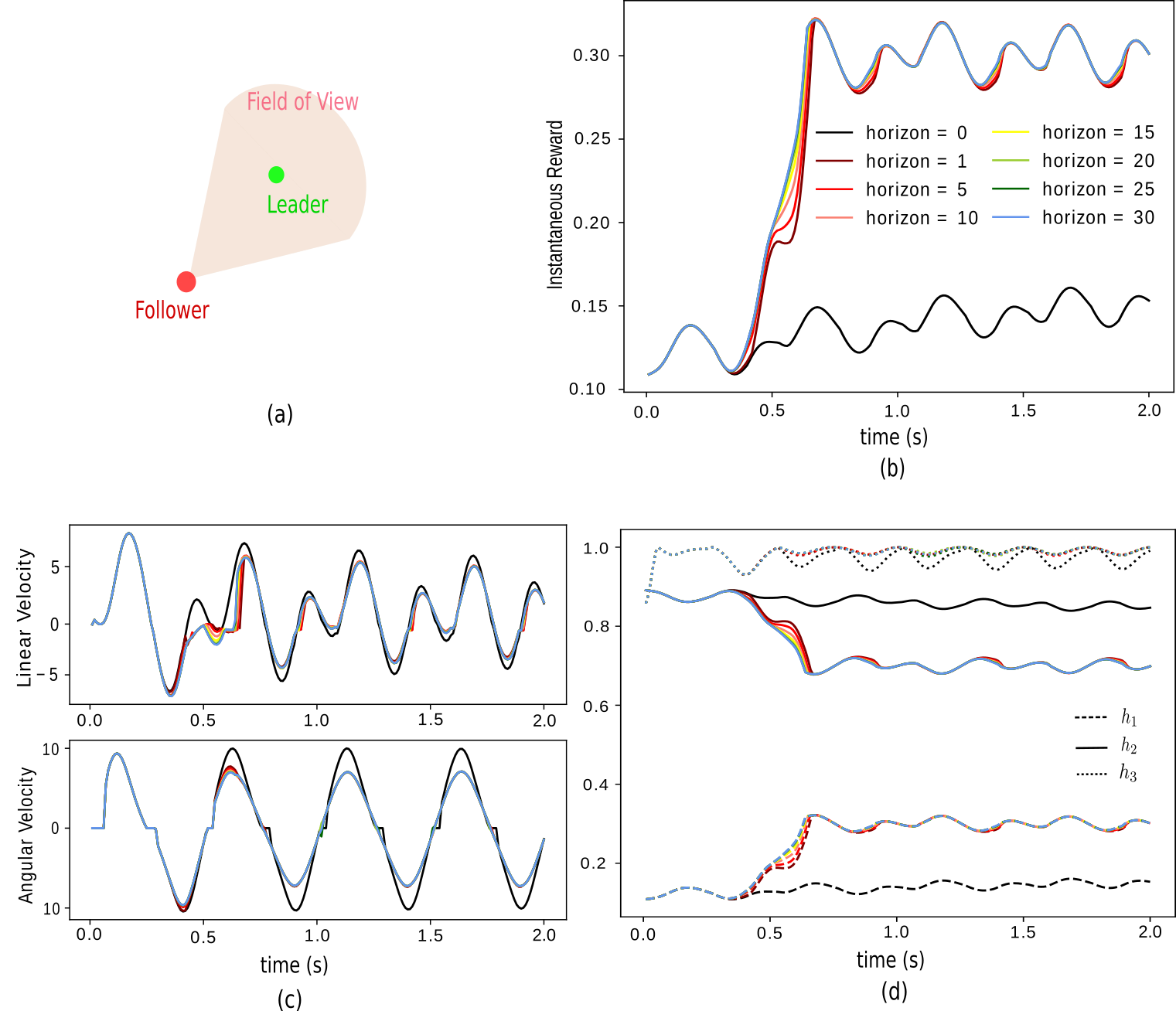}
        \caption{Leader-follower simulation comparing different horizons for high-level optimization. (a) shows the scenario considered and (b) shows the rewards received with time. The adaptive parameter case manages to achieve higher rewards. (c) shows the control input variation and (d) shows variation of barrier functions with time. $h_1$ and $h_2$ represent minimum and maximum distance barriers respectively and converge to the value corresponding to desired relative location of leader.}
    \label{fig:unicycle_plots}
\end{figure}

The leader is made to move with constant horizontal velocity of $\dot{x}=1$ m/s and sinusoidal vertical velocity of $\dot{y} = 12\sin(4\pi t)$ m/s. Parameter values for follower's QP are continuously updated online. At each time instant, a single step of GD is done for a look-ahead horizon that ranges from 1 to 30 in simulations. Fig. \ref{fig:unicycle_plots} shows the difference between the constant parameter case(horizon 0) and the adaptive parameter cases. The rewards received in Eq.(\ref{eq::high_level}) for the proposed algorithm are much higher compared to constant parameter case as seen in Fig.\ref{fig:unicycle_plots}(b). The proposed algorithm improves the system's performance without losing feasibility.

\subsection{Computational Complexity for factored graphs}
A thorough analysis of computational complexity is outside the scope of this paper.
For simplicity, we only evaluate number of gradient computations. For every function dependent on $x_t,u_t$, the number of gradients in backpropagation can be easily deduced to be $4t-1$ by counting edges that need to be traversed when going back to initial state $x_1$ and parameter vector $\theta$. Since we have $N$ hard constraints and one reward at each time $t\in\{1,2..,T\}$, the total number of gradients required are $\sum_{t=1}^{T}(N+1)(4t-1)=(N+1)T(T+1)$. In practice, however, the constraints that are already satisfied by a large margin need not be considered for feasibility maintenance as small changes in parameters will not affect them. Also, the updates may be performed based on linearization about past trajectory instead of predicted one as it would offer huge boost in speed.

\section{Conclusion and Future Work}
This paper details a framework to update parameters of QP controller by employing backpropagation on a computational graph over time. Sensitivity analysis was used for the first time and combined with SQP algorithms to ensure that feasibility of the process over the horizon of interest is maintained or improved while simultaneously achieving performance gains.
Future work will investigate experimental verification on different robotic platforms, convergence rates guarantees, and the effect of uncertain dynamics. Since differentiating through conic programs is not limited to QPs and LPs, more generic optimization based controllers will also be investigated. While the unicycle example showed how algorithm would perform for a nonlinear dynamical system, the algorithm needs to be tested for feasibility on more systems and for infeasibility resulting from control input bounds. We also plan to use more general \classK functions for imposing barrier and lyapunov conditions.



\bibliographystyle{IEEEtrans}
\bibliography{icra.bib}

\begin{thebibliography}{10}
\providecommand{\url}[1]{#1}
\csname url@rmstyle\endcsname
\providecommand{\newblock}{\relax}
\providecommand{\bibinfo}[2]{#2}
\providecommand\BIBentrySTDinterwordspacing{\spaceskip=0pt\relax}
\providecommand\BIBentryALTinterwordstretchfactor{4}
\providecommand\BIBentryALTinterwordspacing{\spaceskip=\fontdimen2\font plus
\BIBentryALTinterwordstretchfactor\fontdimen3\font minus
  \fontdimen4\font\relax}
\providecommand\BIBforeignlanguage[2]{{%
\expandafter\ifx\csname l@#1\endcsname\relax
\typeout{** WARNING: IEEEtran.bst: No hyphenation pattern has been}%
\typeout{** loaded for the language `#1'. Using the pattern for}%
\typeout{** the default language instead.}%
\else
\language=\csname l@#1\endcsname
\fi
#2}}

\bibitem{agrawal2019differentiating}
A.~Agrawal, S.~Barratt, S.~Boyd, E.~Busseti, and W.~M. Moursi,
  ``Differentiating through a cone program,'' \emph{J. Appl. Numer. Optim},
  vol.~1, no.~2, pp. 107--115, 2019.

\bibitem{agrawal2017discrete}
A.~Agrawal and K.~Sreenath, ``Discrete control barrier functions for
  safety-critical control of discrete systems with application to bipedal robot
  navigation.'' in \emph{Robotics: Science and Systems}, vol.~13.\hskip 1em
  plus 0.5em minus 0.4em\relax Cambridge, MA, USA, 2017.

\bibitem{ames2014control}
A.~D. Ames, J.~W. Grizzle, and P.~Tabuada, ``Control barrier function based
  quadratic programs with application to adaptive cruise control,'' in
  \emph{53rd IEEE Conference on Decision and Control}.\hskip 1em plus 0.5em
  minus 0.4em\relax IEEE, 2014, pp. 6271--6278.

\bibitem{ames2016control}
A.~D. Ames, X.~Xu, J.~W. Grizzle, and P.~Tabuada, ``Control barrier function
  based quadratic programs for safety critical systems,'' \emph{IEEE
  Transactions on Automatic Control}, vol.~62, no.~8, pp. 3861--3876, 2016.

\bibitem{amos2017optnet}
B.~Amos and J.~Z. Kolter, ``Optnet: Differentiable optimization as a layer in
  neural networks,'' in \emph{International Conference on Machine
  Learning}.\hskip 1em plus 0.5em minus 0.4em\relax PMLR, 2017, pp. 136--145.

\bibitem{bansal2017hamilton}
S.~Bansal, M.~Chen, S.~Herbert, and C.~J. Tomlin, ``Hamilton-jacobi
  reachability: A brief overview and recent advances,'' in \emph{2017 IEEE 56th
  Annual Conference on Decision and Control (CDC)}.\hskip 1em plus 0.5em minus
  0.4em\relax IEEE, 2017, pp. 2242--2253.

\bibitem{barratt2018differentiability}
S.~Barratt, ``On the differentiability of the solution to convex optimization
  problems,'' \emph{arXiv preprint arXiv:1804.05098}, 2018.

\bibitem{biegler2010nonlinear}
L.~T. Biegler, \emph{Nonlinear programming: concepts, algorithms, and
  applications to chemical processes}.\hskip 1em plus 0.5em minus 0.4em\relax
  SIAM, 2010.

\bibitem{blanchini1999set}
F.~Blanchini, ``Set invariance in control,'' \emph{Automatica}, vol.~35,
  no.~11, pp. 1747--1767, 1999.

\bibitem{boyd2004convex}
S.~Boyd, S.~P. Boyd, and L.~Vandenberghe, \emph{Convex optimization}.\hskip 1em
  plus 0.5em minus 0.4em\relax Cambridge university press, 2004.

\bibitem{breeden2021control}
J.~Breeden, K.~Garg, and D.~Panagou, ``Control barrier functions in
  sampled-data systems,'' \emph{IEEE Control Systems Letters}, 2021.

\bibitem{cohen2020approximate}
M.~H. Cohen and C.~Belta, ``Approximate optimal control for safety-critical
  systems with control barrier functions,'' in \emph{2020 59th IEEE Conference
  on Decision and Control (CDC)}.\hskip 1em plus 0.5em minus 0.4em\relax IEEE,
  2020, pp. 2062--2067.

\bibitem{diamond2016cvxpy}
S.~Diamond and S.~Boyd, ``Cvxpy: A python-embedded modeling language for convex
  optimization,'' \emph{The Journal of Machine Learning Research}, vol.~17,
  no.~1, pp. 2909--2913, 2016.

\bibitem{glotfelter2017nonsmooth}
P.~Glotfelter, J.~Cort{\'e}s, and M.~Egerstedt, ``Nonsmooth barrier functions
  with applications to multi-robot systems,'' \emph{IEEE control systems
  letters}, vol.~1, no.~2, pp. 310--315, 2017.

\bibitem{gros2019towards}
S.~Gros and M.~Zanon, ``Towards safe reinforcement learning using nmpc and
  policy gradients: Part ii-deterministic case,'' \emph{arXiv preprint
  arXiv:1906.04034}, 2019.

\bibitem{panagou2013multi}
D.~Panagou, D.~M. Stipanovi{\v{c}}, and P.~G. Voulgaris, ``Multi-objective
  control for multi-agent systems using lyapunov-like barrier functions,'' in
  \emph{52nd IEEE Conference on Decision and Control}.\hskip 1em plus 0.5em
  minus 0.4em\relax IEEE, 2013, pp. 1478--1483.

\bibitem{stipanovic2012monotone}
D.~M. Stipanovi{\'c}, C.~J. Tomlin, and G.~Leitmann, ``Monotone approximations
  of minimum and maximum functions and multi-objective problems,''
  \emph{Applied Mathematics \& Optimization}, vol.~66, no.~3, pp. 455--473,
  2012.

\bibitem{tits2009feasible}
A.~L. Tits, ``Feasible sequential quadratic programming,'' \emph{Encyclopedia
  of Optimization}, 2009.

\bibitem{xiao2020feasibility}
W.~Xiao, C.~A. Belta, and C.~G. Cassandras, ``Feasibility-guided learning for
  constrained optimal control problems,'' in \emph{2020 59th IEEE Conference on
  Decision and Control (CDC)}.\hskip 1em plus 0.5em minus 0.4em\relax IEEE,
  2020, pp. 1896--1901.

\end{thebibliography}

\end{document}